\documentclass[pdflatex]{article}
\usepackage[fleqn]{amsmath}
\usepackage{ntheorem-hyper}
\usepackage{amssymb}
\usepackage{authblk} 
\usepackage{algorithm}
\usepackage{algorithmicx}
\usepackage{algpseudocode}
\usepackage{graphicx}
\usepackage{graphbox}
\usepackage{subfigure}
\usepackage{svg}
\usepackage{epstopdf}
\usepackage{tabularx}  
\usepackage{booktabs}
\usepackage[numbers, sort]{natbib}
\usepackage{hyperref}
\usepackage[top=2cm, bottom=2cm, left=2cm, right=2cm]{geometry}  
\usepackage{xcolor} 
\usepackage{color}
\hypersetup{urlcolor=blue, citecolor=red}

\usepackage{tikz}
\usetikzlibrary{arrows.meta, positioning, calc}
\usepackage[T1,T2A]{fontenc}
\usepackage[utf8]{inputenc}
\usepackage[russian, german, english]{babel}

\usepackage{authblk}

\newtheorem{theorem}{Theorem}
\newtheorem{prop}{Proposition}
\theorembodyfont{\upshape}
\newtheorem*{proof}{Proof}

\newtheorem{problem}{Problem}

\def\RR{\mathbb{R}}

\def\ZZ{\mathbb{Z}}

\title{On the geometry of Riemannian isometric embeddings}
\author[a]{Dmitri Burago \\ Email: dyb6@psu.edu
	\and 
	Hongda Qiu \\ Email: hjq5042@psu.edu}
\begin{document}
\maketitle
\begin{abstract}
	This note pertains to isometric embeddings endowed with certain geometric properties. We study two embedding problems for a Riemannian manifold $M$ which is diffeomorphic to $\RR^n$ and admits a Bieberbach group $\Gamma$ acting by isometries. The first problem concerns the existence of an isometric embedding of $M$ into a bounded subset of some Euclidean space $\RR^{D_1}$. The second problem seeks a $\Gamma$-equivariant isometric embdding of $M$ into $\RR^{D_2}$. By using a known trick in a novel way, our idea yields results with $D_1 = N+2n$ and $D_2 = N+n$, where $N$ is the Nash dimension of $ M/\Gamma $. Moreover, we also show that an $n$-dimensional smooth manifold, of Nash dimension $N$, can be isometrically embedded into a bounded subset of $\RR^{2N}$.
	
\end{abstract}
\section{Introduction}
\textit{History.} In the early 1950s, John Nash proved that every smooth manifold can be isometrically smoothly embedded into some Euclidean space \cite[Theorem 3, p.~63]{nash1956}. Later, several mathematicians, including Gromov, Rokhlin \cite{Gromov1970}\cite{Gromov1986}, and Günther \cite{Gunther1989}, contributed significantly to deeper understanding of the dimension bounds for the target space. However, not so much is known about the embeddings themselves besides smoothness, and exploring the structure of Nash embeddings remains an intriguing topic. 

We work with $C^\infty$-smoothness unless further specified. It is known that the embedding can be constructed to have the same smoothness as the manifold, and other smoothness classes still contain many problems that seems to be open.



Our objects of interest are Riemannian manifolds with symmetries. Consider an $n$-dimensional Riemannian manifold $(M,g)$. The symmetries of $M$ are represented by a group $\Gamma$ acting on $M$ by isometries. A basic example is the universal cover of a Riemannian torus. In this case, the manifold is acted upon by $\ZZ^n$. In our work, we consider a more general case when $(M,g)$ is diffeomorphic to $\RR^n$ and admits a group $\Gamma$ acting properly discontinuously. We also assume that $\Gamma$ acts co-compactly. Such groups are called Bieberbach groups. 

In this work, we separately study two aspects for isometric embeddings of $M$: boundedness and equivariance, and construct embeddings having these properties. In our results, the dimensions of target spaces are in terms of $n$ and the \textit{Nash dimension} of $M$, which is defined as the smallest value $N = N(M)$ such that $ M $ equipped with any Riemannian metric can be isometrically embedded into $\RR^{N}$.


\textit{Boundedness.} We study the following problem:
\begin{problem}[Isometric embeddings into a bounded domain]
	Can $(M,g)$ be isometrically embedded into a bounded subset of some Euclidean space $\RR^{D_1}$ for some good enough $D_1$ independent of the metric $g$?
\end{problem}  
The existence of bounded embeddings and immersions has been studied historically. For instance, Rozendorn \cite[p.~117]{rozendorn1992} proved Hadamard's conjecture that there are no bounded complete surfaces of negative curvature in $\RR^3$. Based on this result, he further showed that there are no $C^2$-smooth closed surfaces with a metric of non-positive curvature. A negative result was given by one of the author \cite{burago1984unboundednessofhorn}, stating that the $C^2$-immersion of tubes having finite total intrinsic curvature in $\RR^3$ cannot be bounded. More results are surveyed in \cite{Burago2024} and \cite{rozendorn1992}.

\textit{Equivariance.} Recall that $\Gamma$ acts on $M$ by isometries. An embedding $M\to\RR^D$ is said to be \textit{$\Gamma$-equivariant} if $\Gamma$ acts not only on the image of the embedding but also on the entire target space $\RR^D$ as in the commutative diagram in Figure~\ref{figure equivariant-map}. Our key questions are: 
\begin{problem}[Equivariant isometric embeddings]
	\label{problem equivariant}
	Can $(M,g)$ be equivariantly isometrically embedded into some Euclidean space $\RR^{D_2}$ for some good enough $D_2$ independent of the Riemannian metric on $M$? 
\end{problem}

\begin{figure}[htbp]
	\centering
	\begin{tikzpicture}[
		node distance=2.5cm, auto,
		>=Latex,
		every node/.style={font=\small}
		]
		\node (M)   {$M$};
		\node (M')  [right=of M] {$M$};
		\node (N)   [below=of M] {$\RR^D$};
		\node (N')  [below=of M'] {$\RR^D$};
		
		\draw[->] (M) -- node[left] {$u$} (N);
		\draw[->] (M') -- node[right] {$u$} (N');
		\draw[->] (M) -- node[above] {$ g\in \Gamma $} (M');
		\draw[->] (N) -- node[below] {$ g\in \Gamma $} (N');
		
		\node at ($(M)!0.5!(N')$) [align=center] {\\[3pt]$u\circ g = g\circ u$};
	\end{tikzpicture}
	\caption{A diagram illustrating a $\Gamma$-equivariant embedding $ u:M\to \RR^D $. Here we use the same notation $g$ for an element of $\Gamma$ and its corresponding actions.}
	\label{figure equivariant-map}
\end{figure}

The term ``equivariant isometric embedding" was introduced in a number of works. Some progress has been made in special cases. Regarding the existence problem, there are many positive results based on Mostow's theorem \cite[Theorem 6.1]{mostwo1956} that a compact manifold acted upon by a compact Lie group can be equivariantly embedded into some Euclidean space. For instance, Moore \cite{Moore1976} proved the existence of equivariant isometric embedding when $M$ is a compact Riemannian homogeneous space and $\Gamma$ is its isometry group (for sufficiently large $D$). Later, in a joint work with Schlafly \cite{Moore1980}, he proved a similar theorem for compact manifolds acted upon by compact Lie groups. In the same work, they also established results for embeddings into other type of target spaces. 

For specific examples, these works illustrate how Nash's original approach can be modified to impose additional conditions on the isometric embedding beyond smoothness. They shed some light on the construction of general equivariant isometric embeddings determined by the topology of $M$. This is in the focus of our study. 




The discussion on equivariance can be divided into several paradigms as follows.
\begin{enumerate}
	\item The group $\Gamma$ just cannot act on Euclidean spaces. For instance, a discrete group of exponential growth.
	
	
	\item The group $\Gamma$ is a compact Lie group and can act on Euclidean spaces. Then the answer is covered by Moore and Schlafly's results,  though the bound for $D$ depends on both $M$ and $\Gamma$, stating that such an equivariant isometric embedding exists though there is no universal dimension bound for the target space \cite[Section 6]{Moore1980}. There are further results on specific examples such as compact symmetric spaces \cite{Eschenburg2022} and Hermitian symmetric spaces \cite{clozel2007}, where the bound for $D$ depends only on $M$. 
	
	\item The group $\Gamma$ acts properly discontinuously on Euclidean spaces. Little is explored about this case, and the precise bounds for $D$ are unknown. This is the direction pursued in our current work.
\end{enumerate}

\textit{Results.} Recall that $N$ is the Nash dimension of the quotient manifold $M/\Gamma$. In this work, we show that it is possible to take $D_1 = N+2n$ and $D_2 = N+n$. In other words, the dimensions of target spaces depend only on $M$.

\begin{theorem}[Isometric embedding into a bounded region]
	\label{theorem bounded}
	Let $M$ be a Riemannian manifold that is diffeomorphic to $\RR^n$ and admits an action of a Bieberbach group $\Gamma$. Then there exists an isometric embedding ${u}: (M,g)\to\RR^{N+2n}$ such that ${u}(M)$ is a bounded subset of $\RR^{N+2n}$.
\end{theorem}
Along our discussion, we also prove the following theorem for general Riemannian manifolds:
\begin{theorem}
	\label{theorem general bounded}
	Let $M$ be a Riemannian manifold and $N$ its Nash dimension. Then there exists an isometric embedding of $M$ into a bounded region in $\RR^{2N}$.
\end{theorem}
We remark that one might probably decrease this dimension bound from $2N$ to $N+1$, though we did not check the details yet.
\begin{theorem}[Equivariant isometric embedding]
	\label{theorem equivariant}
	Let $M$ be a Riemannian manifold that is diffeomorphic to $\RR^n$ and admits an action of a Bieberbach group $\Gamma$. Then there exists a $\Gamma$-equivariant isometric embedding ${{u}}:(M,g)\to\RR^{N+n}$.
\end{theorem}
\textit{Remarks.} Our results are still far from satisfactory, since the Nash embedding theorem is still used as a ``blackbox" to guarantee the existence of some auxiliary embedding in our constructions of embeddings. It remains unclear what the optimal values of $D_1,D_2$ are. Our humble hope is to improve the bounds to $N$ by adjusting Nash and Günther's techniques to our situation.


Smoothness is yet another sensitive aspect for further study. In this work, we stay consistent with historical results in $C^\infty$. However, one possibly can obtain a very different story in the $C^2$ case regarding boundedness and equivariance, and the $C^k,k\geq3$ case remains obscure.






\section{Proofs}
We aim to construct isometric embeddings as mentioned in Problems 1 and 2. Both constructions begin in a similar manner as follows.

For convenience, in what remains, we use the same notations for metrics on $M/\Gamma$ and their lifts to $M$.

We start by decomposing $g$ into the sum of a Riemannian metric ${g_1}$ and a small flat metric ${g_2}$ on $M$ (hence $(M,g_2)$ is isometric to $\RR^n$). Since $M/\Gamma$ is compact, we can take ${g_2}$ small enough so that $ {g_1} = g-{g_2} $ remains positive definite. 


By the construction above, both ${g_1}$ and ${g_2}$ are smooth and define Riemannian metrics on $M$. It follows that $ (M,g) $ is isometric to an $n$-dimensional submanifold of the product manifold $ (M,{g_1})\times(M,{g_2}) $ via the diagonal map $x\mapsto(x,x)$. 

We construct our embeddings as follows. We first immerse $ (M,{g_1}) $ as a covering of the quotient $(M/\Gamma,{g_1})$ which is itself isometrically embedded into some Euclidean space via Nash embedding. Next, we isometrically embed $(M,{g_2})$ into a Euclidean space. To be specific, it is embedded into a bounded subset of $\RR^{2n}$ for Problem 1 and $\RR^n$ for Problem 2. We then pair the immersion of $(M,{g_1})$ and embedding of $(M,{g_2})$ to obtain an isometric embedding, and show that these embeddings resolve Problems 1 and 2.

Let us formalize the discussion above in the following two propositions. 
\begin{prop}
	\label{proposition q1}
	There exists a $\Gamma$-equivariant isometric immersion ${u_1}:(M,{g_1})\to\RR^N$.
\end{prop}
\begin{proof}
	This is obvious. We can take the composition of the covering map $M\to M/\Gamma$ and the Nash embedding $M/\Gamma \to\RR^N$; that is, $ {u_1}: (M,{g_1})\to (M/\Gamma,g_1)\to \RR^N$. 
\end{proof}

\begin{prop}
	\label{proposition q2}
	There exists an isometric embedding ${u_2}:\RR^n\to\RR^{2n}$ such that the image ${u_2}(\RR^n) $ is bounded.
\end{prop}
\begin{proof}
	
	The real line $\RR^1$ can be isometrically embedded into a bounded subset of $\RR^2$ via a map $\rho:\RR^1\to\RR^2$, for instance, as a spiral curve confined in an annulus (Figure~\ref{figure spiral curve}). 
	\begin{figure}[htbp]
		\centering
		\caption{A spiral curve winding in an annulus}
		\label{figure spiral curve}
		\includegraphics[width=0.75\linewidth]{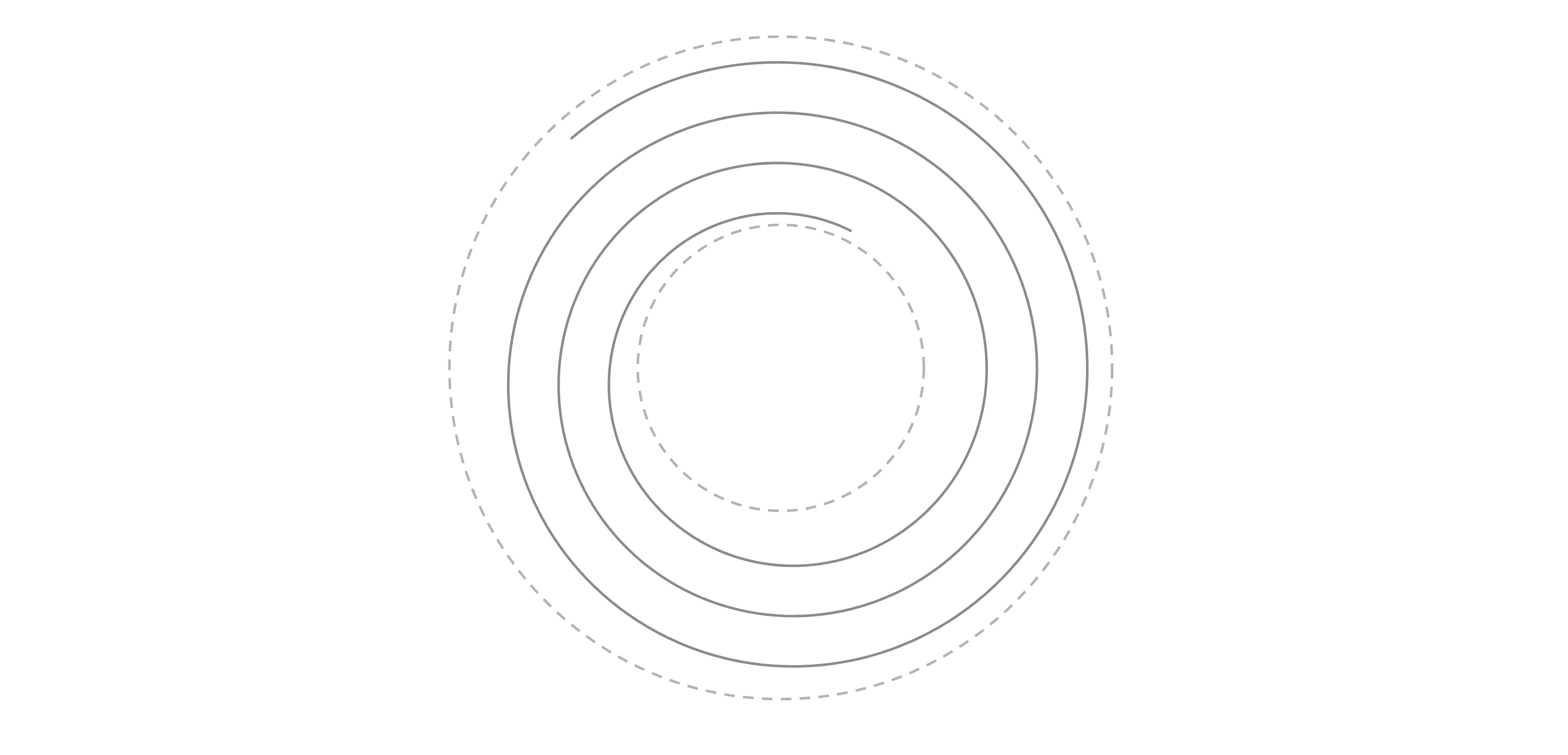}
	\end{figure}
	
	
	We can define ${u_2}:\RR^n\to\RR^{2n}$ as follows:
	\begin{equation}
		x\mapsto (\rho(x_1),...,\rho(x_n)).
	\end{equation}
	This map isometrically embeds $\RR^n$ into a bounded subset of $\RR^{2n}$. 
\end{proof}

\begin{proof}[Theorem~\ref{theorem general bounded}]
	We start with the Nash embedding $M\to\RR^{N}$, and compose it with a bounded isometric embedding $u_2:\RR^{N}\to\RR^{2N}$ following Proposition~\ref{proposition q2}. 
	
\end{proof}

\begin{proof}[Theorem~\ref{theorem bounded}]
	With $ {u_1}:(M,{g_1})\to \RR^N $ from Proposition~\ref{proposition q1} and $ {u_2}:\RR^n\to \RR^{2n} $ from Proposition~\ref{proposition q2}, we can isometrically embed $(M,g)$ into $ (M,{g_1})\times(M,{g_2}) $ via the diagonal map $ i:(M,g) \to (M,{g_1})\times(M,{g_2}), x\mapsto(x,x) $. Since $(M,g_2)$ is isometric to $\RR^n$, we can define
	\begin{equation}
		{u} = ({u_1},{u_2})\circ i: (M,g) \to \RR^{N}\times\RR^{2n}, x\mapsto ({u_1}(x),{u_2}(x)).
	\end{equation}
	Since ${u_2}$ is an embedding (hence injective), $ {u} $ is an injective immersion. To show that ${u}$ is an embedding, it suffices to show that it is proper, i.e, the preimage of every compact subset is compact. 
	\footnote{In general, an injective immersion need not to be an embedding. The correct statement is that a proper injective immersion is exactly an embedding \cite[Proposition 4.22, p.~87]{lee2012smoothmanifolds}. For an example that an injective immersion fails to be an embedding, see Example 4.19 (p.~86) of the same reference.}
	By our construction, for any compact set $K$ in ${u}(M)$, its preimage is completely determined by its projection $K'$ onto the $\RR^{2n}$ factor; specifically, $ {u}^{-1}(K) = {u_2}^{-1}(K') $. Since $K'$ is compact, it follows that ${u}^{-1}(K)$ is also compact.
	
	Since both ${u_1},{u_2}$ are Riemannian isometries onto their images, $E$ is an isometric embedding. Therefore, the image ${u}(M) \subset {u_1}(M)\times{u_2}(M)$ is bounded.
\end{proof}
Next, we turn to the construction of an equivariant isometric embedding of $M$. We achieve this by pairing the isometric immersion ${u_1}:(M,{g_1})\to\RR^N$ from Proposition $\ref{proposition q1}$ with the isometry between $(M,{g_2})$ and $\RR^n$.


\begin{proof}[Theorem~\ref{theorem equivariant}]
	Note that $(M,{g_2})$ is isometric to $\RR^n$. We define 
	\begin{equation}
		{{u}} = ({u_1},id)\circ i:M\to \RR^{N}\times\RR^n, x\mapsto ({u_1}(x),x),
	\end{equation}
	where $i$ is the diagonal map $ i:(M,g) \to (M,{g_1})\times(M,{g_2}), x\mapsto(x,x) $. Since the map $x\mapsto x$ is injective, ${{u}}$ is an injective immersion. Following a step similar to the one at the beginning of the proof of Theorem~\ref{theorem bounded}, ${{u}}$ is an isometric embedding. Finally, since $u_1$ is $\Gamma$-equivariant, so is $u$. 
	
%
%
\end{proof}

\newpage
\bibliographystyle{unsrt}
\bibliography{ref}
\end{document}